\documentclass{amsart}

\usepackage{amsmath,amsthm,amsfonts,amscd,amssymb}
\usepackage{verbatim}
\usepackage{color} 
\usepackage{graphicx}
\usepackage{float}

\usepackage[square,sort,comma,numbers]{natbib}

\newtheorem{thm}{Theorem}[section]

\newtheorem{lem}[thm]{Lemma}
\newtheorem{prop}[thm]{Proposition}

\newtheorem{quest}{Question}

\newtheorem{ex}{Example}

\theoremstyle{definition}
\newtheorem{defn}{Definition}[section]

\theoremstyle{remark}
\newtheorem{rem}{Remark}[section]

\begin{document}

\title[Canonical basis twists of ideal lattices]{Canonical basis twists of ideal lattices from real quadratic number fields}

\author{Mohamed Taoufiq Damir and Lenny Fukshansky}\thanks{Damir was partially supported by the Academy of Finland through grants \#276031, \#282938, and \#303819 awarded to C. Hollanti. Fukshansky was partially supported by the Simons Foundation grant \#519058.}

\address{Department of Mathematics and Systems Analysis, Aalto University, P.O. Box 11100, FI-00076 Aalto, Finland}
\email{mohamed.damir@aalto.fi}

\address{Department of Mathematics, 850 Columbia Avenue, Claremont McKenna College, Claremont, CA 91711, USA}
\email{lenny@cmc.edu}


\begin{abstract}
Ideal lattices in the plane coming from real quadratic number fields have been investigated by several authors in the recent years. In particular, it has been proved that every such ideal has a basis that can be twisted by the action of the diagonal group into a Minkowski reduced basis for a well-rounded lattice. We explicitly study such twists on the canonical bases of ideals, which are especially important in arithmetic theory of quadratic number fields and binary quadratic forms. Specifically, we prove that every fixed real quadratic field has only finitely many ideals whose canonical basis can be twisted into a well-rounded or a stable lattice in the plane. We demonstrate some explicit examples of such twists. We also briefly discuss the relation between stable and well-rounded twists of arbitrary ideal bases.
\end{abstract}

\maketitle

\def\A{{\mathcal A}}
\def\AA{{\mathfrak A}}
\def\B{{\mathcal B}}
\def\C{{\mathcal C}}
\def\D{{\mathcal D}}
\def\E{{\mathcal E}}
\def\F{{\mathcal F}}
\def\Ff{{\mathfrak F}}
\def\G{{\mathcal G}}
\def\x{{\mathcal H}}
\def\I{{\mathcal I}}
\def\J{{\mathcal J}}
\def\K{{\mathcal K}}
\def\kk{{\mathfrak K}}
\def\L{{\mathcal L}}
\def\LL{{\mathfrak L}}
\def\M{{\mathcal M}}
\def\O{{\mathcal O}}
\def\W{{\mathcal W}}
\def\CC{{\mathfrak C}}
\def\mm{{\mathfrak m}}
\def\MM{{\mathfrak M}}
\def\OO{{\mathfrak O}}
\def\P{{\mathcal P}}
\def\R{{\mathcal R}}
\def\s{{\mathcal S}}
\def\V{{\mathcal V}}
\def\X{{\mathcal X}}
\def\XX{{\mathfrak X}}
\def\Y{{\mathcal Y}}
\def\Z{{\mathcal Z}}
\def\H{{\mathcal H}}
\def\cee{{\mathbb C}}
\def\pee{{\mathbb P}}
\def\que{{\mathbb Q}}
\def\real{{\mathbb R}}
\def\zed{{\mathbb Z}}
\def\hyp{{\mathbb H}}
\def\aaa{{\mathbb A}}
\def\Nn{{\mathbb N}}
\def\ff{{\mathbb F}}
\def\kk{{\mathfrak K}}
\def\qbar{{\overline{\mathbb Q}}}
\def\kbar{{\overline{K}}}
\def\ybar{{\overline{Y}}}
\def\kkbar{{\overline{\mathfrak K}}}
\def\ubar{{\overline{U}}}
\def\eps{{\varepsilon}}
\def\ahat{{\hat \alpha}}
\def\bhat{{\hat \beta}}
\def\gt{{\tilde \gamma}}
\def\h{{\tfrac12}}
\def\dd{{\partial}}
\def\baa{{\boldsymbol \alpha}}
\def\bfa{{\boldsymbol a}}
\def\bfb{{\boldsymbol b}}
\def\be{{\boldsymbol e}}
\def\bei{{\boldsymbol e_i}}
\def\bff{{\boldsymbol f}}
\def\bc{{\boldsymbol c}}
\def\bm{{\boldsymbol m}}
\def\bk{{\boldsymbol k}}
\def\bi{{\boldsymbol i}}
\def\bl{{\boldsymbol l}}
\def\bq{{\boldsymbol q}}
\def\bu{{\boldsymbol u}}
\def\bt{{\boldsymbol t}}
\def\bs{{\boldsymbol s}}
\def\bfu{{\boldsymbol u}}
\def\bv{{\boldsymbol v}}
\def\bw{{\boldsymbol w}}
\def\bx{{\boldsymbol x}}
\def\bX{{\boldsymbol X}}
\def\bz{{\boldsymbol z}}
\def\bwy{{\boldsymbol y}}
\def\bY{{\boldsymbol Y}}
\def\bL{{\boldsymbol L}}
\def\ba{{\boldsymbol a}}
\def\bb{{\boldsymbol b}}
\def\bet{{\boldsymbol\eta}}
\def\bxi{{\boldsymbol\xi}}
\def\bo{{\boldsymbol 0}}
\def\bol{{\boldkey 1}_L}
\def\ep{\varepsilon}
\def\p{\boldsymbol\varphi}
\def\q{\boldsymbol\psi}
\def\rank{\operatorname{rank}}
\def\aut{\operatorname{Aut}}
\def\lcm{\operatorname{lcm}}
\def\sgn{\operatorname{sgn}}
\def\spn{\operatorname{span}}
\def\md{\operatorname{mod}}
\def\Norm{\operatorname{Norm}}
\def\dim{\operatorname{dim}}
\def\det{\operatorname{det}}
\def\Vol{\operatorname{Vol}}
\def\rk{\operatorname{rk}}
\def\ord{\operatorname{ord}}
\def\ker{\operatorname{ker}}
\def\div{\operatorname{div}}
\def\Gal{\operatorname{Gal}}
\def\GL{\operatorname{GL}}
\def\SL{\operatorname{SL}}
\def\p{\operatorname{p}}
\def\q{\operatorname{q}}
\def\t{\operatorname{t}}
\def\hs{{\hat \sigma}}
\def\chr{\operatorname{char}}
\def\diag{\operatorname{diag}}

\section{Introduction}
\label{intro}

Euclidean lattices, metrized free $\zed$-modules in real vector spaces, are central to number theory and discrete geometry. In addition to their arithmetic and geometric appeal, lattices are also key to discrete optimization, often providing solutions to classical optimization questions like sphere packing, covering and kissing number problems. They are also extensively used in applied areas, such as coding theory, cryptography and other areas of digital communications; see the famous book by Conway and Sloane~\cite{conway} for a wealth of information on the theory of lattices and their numerous connections and applications. In this paper, we focus on certain geometric properties of lattices that are of utmost interest for both, theory and applications. One of the main sources of lattices possessing these special properties is the classical Minkowski construction from ideals in rings of integers of algebraic number fields: not only does this construction allow for a nice and compact description of the resulting lattices, but also algebraic properties of an ideal often inform the geometry of the corresponding lattice.

Let $L \subset \real^n$ be a free $\zed$-module of rank $k \leq n$, then 
$$L = B\zed^n = \spn_{\zed} \{ \bb_1,\dots,\bb_k \},$$
where $\bb_1,\dots,\bb_n \in \real^n$ are $\real$-linearly independent basis vectors for $L$ and $B = (\bb_1\ \dots\ \bb_n) \in \GL_n(\real)$ is the corresponding basis matrix. Then for any $U \in \GL_n(\zed)$, $BU$ is also a basis matrix for $L$, i.e. $L = BU\zed^n$ for any $U \in \GL_n(\zed)$. 

Let $f(\bx) = f(x_1,\dots,x_n) \in \real[x_1,\dots,x_n]$ be a positive definite quadratic form, i.e. $f(\bx) = \bx^t Q \bx$, where $Q$ is an $n \times n$ symmetric positive definite coefficient matrix for $f$. Then $Q = A^t A$ for some $A \in \GL_n(\real)$, so
\begin{equation}
\label{qf}
f(\bx) = (A\bx)^t (A\bx).
\end{equation}
Thus we will use notation $f_A$ for the form $f$ as in~\eqref{qf} with coefficient matrix $A^t A$ for some $A \in \GL_n(\real)$. 

We will use the term {\it lattice} to refer to a pair $(L,f_A)$, where the form $f_A$ is used to define the norm of vectors in $L$. If $A=I_n$, the $n \times n$ identity matrix, then $f_A$ is $\|\ \|^2$, the square of the usual Euclidean norm, in which case we simply write $L$ instead of $(L, f_{I_n})$. In general, let $\bx = B \bwy \in L$, i.e. $\bwy \in \zed^n$, so
\begin{equation}
\label{qA}
f_A(\bx) = (AB\bwy)^t (AB\bwy) = \bwy^t (AB)^t (AB) \bwy = \| (AB) \bwy \|^2.
\end{equation}
In other words, a lattice $(L,f_A)$ can be identified with the lattice~$AL$, which we will refer to as the {\it twist} of $L$ by $A$. Now, for a lattice $L$ we define its rank, $\rk(L)$, to be the cardinality of its basis, and its determinant
$$\det(L) := \left| \det (B^{\top} B) \right|^{1/2},$$
where $B$ is a basis matrix for~$L$. 

There are two important classes of lattices we will discuss. A lattice $L$ is called {\it well-rounded} (abbreviated WR) if there exist $n$ linearly independent vectors $\bc_1,\dots,\bc_n \in L$ such that
$$\|\bc_1\| = \dots = \|\bc_n\| = |L|,$$
where $|L| := \min_{\bx \in L \setminus \{\bo\}} \|\bx\|.$ On the other hand, $L$ is called {\it stable} if for each sublattice $L' \subseteq L$, 
$$\det(L)^{1/rk(L)} \leq \det(L')^{1/rk(L')},$$
and $L$ is called unstable otherwise. Well-roundedness and stability are independent properties for lattices of rank greater than two: WR lattices can be unstable and stable lattices do not have to be WR. On the other hand, in the plane WR lattices form a proper subset of stable lattices~\cite{lenny:stable}. There is an important equivalence relation on the space of lattices: two lattices are called {\it similar} if they are related by a dilation and an orthogonal transformation; geometric properties, like well-roundedness and stability are preserved under the similarity, hence we can talk about WR or stable similarity classes of lattices.

Both of these types of lattices are very important in reduction theory of algebraic groups. In particular, they were studied in the context of the diagonal group action on the space of lattices. Let
\begin{equation}
\label{A}
\A = \left\{ A = (a_{ij}) \in \GL_n(\real) : a_{ij} = 0\ \forall\ i \neq j,\ a_{ii} > 0\ \forall\ i,\ \prod_{i=1}^n a_{ii} = 1 \right\}
\end{equation}
be the group of real positive diagonal matrices with determinant 1. This group acts on the space of lattices in~$\real^n$ by left multiplication: $L \mapsto AL$ for each $A \in \A$ and lattice $L \subset \real^n$. A celebrated result of McMullen~\cite{mcmullen} in connection with his work on Minkowski's conjecture asserts that any bounded $\A$-orbit of lattices contains a well-rounded lattice. Inspired by McMullen's work, Shapira and Weiss proved~\cite{weiss} that the orbit closure of a lattice under the action of $\A$ also contains a stable lattice. 

Throughout this note, let $K$ be a totally real number field of degree $n$ with the ring of integers $\O_K$, $\sigma_1,\dots,\sigma_n : K \hookrightarrow \real$ be the embeddings of $K$, and define the Minkowski embedding 
$$\sigma_K = (\sigma_1,\dots,\sigma_n) : K \hookrightarrow \real^n.$$
Let $I \subset K$ be a full module, i.e. a $\zed$-module of rank $n$, for instance $\O_K$ or any (fractional) ideal. We define $L_K(I) := \sigma_K(I)$, i.e. the image of $I$ under Minkowski embedding, then $AL_K(I)$ is a lattice in $\real^n$ for any $A \in \GL_n(\real)$. 

\begin{prop} \label{WR2} With notation as above, there exist $A, B \in \A$ such that the lattice $AL_K(I)$ is WR and $BL_K(I)$ is stable.
\end{prop}

\proof
By the similarity relation, we can assume that $L_K(I)$ is unimodular without loss of generality. Consider the orbit of $L_K(I)$ under the action of $\A$, i.e. the set
$$\A L_K(I) = \left\{ A L_K(I) : A \in \A \right\}.$$
Since $L_K(I)$ comes from a full module in a totally real number field, Theorem~3.1 of~\cite{mcmullen} implies that the orbit $\A L_K(I)$ is compact. Then Theorem~1.3 of~\cite{mcmullen} implies that this orbit contains a WR lattice and Theorem~1.1 of~\cite{weiss} implies that it also contains a stable lattice, i.e. there exist some $A,B \in \A$ such that the lattice $AL_K(I)$ is WR and the lattice $BL_K(I)$ is stable.
\endproof

\begin{rem} The proofs of Theorem~1.3 of~\cite{mcmullen} and of Theorem~1.1 of~\cite{weiss} are not constructive, meaning that they do not help to explicitly find $A, B \in \A$ such that the lattice $AL_K(I)$ is WR and $BL_K(I)$ is stable.
\end{rem}
\medskip

Now let $I \subseteq \O_K$ be an ideal. Let $\alpha \in K$ be totally positive, i.e. $\sigma_i(\alpha) > 0$ for all $1 \leq i \leq n$, then the pair $(I,\alpha)$ gives rise to the {\it ideal lattice} $(L_K(I),f_{A(\alpha)})$, where
$$A(\alpha) = \begin{pmatrix} \sqrt{\sigma_1(\alpha)} & \hdots & 0 \\ \vdots & \ddots & \vdots \\ 0 & \hdots & \sqrt{\sigma_n(\alpha)} \end{pmatrix}.$$
Ideal lattices have been extensively studied by a number of authors, especially E. Bayer-Fluckiger (see~\cite{bayer2} for a survey of this topic); a systematic study of WR ideal lattices has been initiated in~\cite{kate-me} and~\cite{fletcher-jones}, where it was shown that WR ideal lattices $(L_K(I),\|\ \|^2)$ from quadratic number fields are relatively sparse. The situation is different if we allow a more general quadratic norm form $f_A$.  Notice that $A(\alpha) \in \A$ if and only if $\alpha \in \O_K$ is a totally positive unit. More generally, $\A_1(K) := \left\{ A(\alpha) : \alpha \in K \text{ is totally positive} \right\}$ is also a multiplicative group of (positive) diagonal matrices so that
$$\A \cap \A_1(K) = \left\{ A(\alpha) : \alpha \in K \text{ is totally positive and } \Nn_K(\alpha) = 1 \right\},$$
where $\Nn_K$ stands for the number field norm on $K$. In fact, notice that the set
$$\A'_1(K) = \left\{ \Nn_K(\alpha)^{-1/n} \A(\alpha) : A(\alpha) \in \A_1(K) \right\}$$
is a proper subset of $\A$, while the lattices $(L_K(I),f_{A(\alpha)})$ and $(L_K(I),f_{\Nn_K(\alpha)^{-1/n} A(\alpha)})$ are scalar multiples of each other, and hence have the same properties. 

\begin{quest}\label{question1} By Proposition~\ref{WR2}, we know that, given $L_K(I)$, there exists $A, B \in \A$ such that $AL_K(I)$ is WR and $BL_K(I)$ is stable, but do there exist such $A, B \in \A_1(K)$?
\end{quest}

In this note, we consider the two-dimensional situation, for which this question was answered in the affirmative for WR lattices in~\cite{bayer-nebe} (specifically, see Corollary~3.1). Further, the authors construct explicit examples of principal ideal lattices (i.e., those coming from full rings of integers) similar to the square and the hexagonal lattices. Of course, this immediately implies the same affirmative answer for stable lattices, since in two dimensions WR lattices are stable; in fact, unlike the WR situation, there are infinitely many stable ideal lattices from real quadratic number fields even without twisting by a matrix, as proved in~\cite{lenny:stable}. Our goal here is to further explicitly investigate well-roundedness and stability properties of planar ideal lattices. First we need some more notation. 

\begin{defn} Given a particular basis matrix $B$ for a lattice $L$, we will say that $B$ is {\it WR twistable} (respectively, {\it stable twistable}) if there exists $A \in \A$ such that $A B \zed^n$ is WR (respectively, stable) with $AB$ being the shortest basis matrix, as discussed in Section~\ref{setup}. We will refer to the resulting lattices as a {\it WR twist} (respectively, {\it stable twist}) of the original lattice by the matrix $A$, respectively. 
\end{defn}

While WR twistable bases have been discussed in~\cite{taoufiq}, to the best of our knowledge stable twistable bases have not previously been investigated. A criterion to determine whether a given basis for an ideal lattice from a real quadratic number field $K$ is WR twistable or not was given in~\cite{taoufiq}: this criterion in particular implies there can be only finitely many such bases for a fixed ideal lattice. Further, Proposition~7 of~\cite{taoufiq} asserts that a basis is WR twistable by some $A \in \A$ if and only if it is WR twistable by some $A(\alpha) \in \A_1(K)$, and an explicit construction of such corresponding twists is given. 

In this note, we focus on the twistable properties of the particularly important kind of basis for ideals in quadratic number fields, the so-called canonical basis. We always talk about twists by matrices $A(\alpha) \in \A_1(K)$ for some totally positive $\alpha \in K$, referring to such twists as a twist by~$\alpha$.

Let $D \in \zed_{>0}$ be squarefree and let $K = \que(\sqrt{D})$. Let $I \subseteq \O_K$ be an ideal. Notice that $\O_K=\zed[\delta]$, where
\begin{equation}
\label{delta}
\delta = \left\{ \begin{array}{ll}
- \sqrt{D} & \mbox{if } D \not\equiv 1\ (\md 4) \\
\frac{1-\sqrt{D}}{2} & \mbox{if } D \equiv 1\ (\md 4) \\
\end{array}
\right.
\end{equation}
The embeddings $\sigma_1, \sigma_2 : K \to \real$ are given by
$$\sigma_1(x+y\sqrt{D}) = x+y\sqrt{D},\ \sigma_2(x+y\sqrt{D}) = x-y\sqrt{D}$$
for each $x+y\sqrt{D} \in K$. The number field norm on $K$ is given by
$$\Nn_K(x+y\sqrt{D}) = \sigma_1(x+y\sqrt{D}) \sigma_2(x+y\sqrt{D}) = \left( x+y\sqrt{D} \right) \left( x-y\sqrt{D} \right).$$
Now $I \subseteq \O_K$ is an ideal if and only if 
\begin{equation}
\label{I_abg}
I = \{ ax + (b+g\delta)y : x,y \in \zed \},
\end{equation}
for some $a,b,g \in \zed_{\geq 0}$ such that
\begin{equation}
\label{abg}
b < a,\ g \mid a,b,\text{ and } ag \mid \Nn(b+g\delta).
\end{equation}
Such integral basis $a,b+g\delta$ is unique for each ideal $I$ and is called the {\it canonical basis} for $I$ (see Section~6.3 of~\cite{buell} for a detailed exposition): it is important in the arithmetic theory of binary quadratic forms and quadratic number fields. For instance, canonical basis is used to determine reduced ideals and compute the ideal class group in quadratic fields (see, e.g., the classical work of H. H. Mitchell~\cite{mitchell} as well as a more recent paper~\cite{yamamoto}), and it is employed for number field computations in several modern computer algebra systems; it has also been used in the previous study of geometric properties of ideal lattices in the plane (see \cite{kate-me}, \cite{fletcher-jones}, \cite{lenny:stable}, \cite{taoufiq}). It is then easy to check that $L_K(I) = B \zed^2$, where 
\begin{equation}
\label{B1}
B = \begin{pmatrix} a & b-g\sqrt{D} \\ a & b+g\sqrt{D} \end{pmatrix},
\end{equation}
when $D \not\equiv 1 (\md 4)$, and
\begin{equation}
\label{B2}
B = \begin{pmatrix} a & \frac{2b+g}{2} - \frac{g\sqrt{D}}{2} \\ a &\frac{2b+g}{2} + \frac{g\sqrt{D}}{2} \end{pmatrix},
\end{equation}
when $D \equiv 1 (\md 4)$. Then any other basis matrix for $L_K(I)$ is of the form $BU$ for some $U \in \GL_2(\zed)$. We can now state our result about WR $K$-twists of canonical bases of planar ideal lattices from real quadratic number fields.
\medskip

\begin{thm} \label{main_wr} Let $K = \que(\sqrt{D})$ be a real quadratic number field and $\O_K$ its ring of integers. There can be at most finitely many ideals in $\O_K$, up to similarity of the resulting ideal lattices, with the canonical basis being WR twistable: if this is the case for some ideal $I$ with canonical basis $a, b+g\delta$, then
$$b < \left\{ \begin{array}{ll}
g\sqrt{D} & \mbox{if $D \not\equiv 1\ (\md 4)$,} \\
\frac{(\sqrt{D}-1)g}{2} & \mbox{if $D \equiv 1\ (\md 4)$.}
\end{array}
\right.$$
The canonical basis for $\O_K$ is WR twistable if and only if $K=\que(\sqrt{5})$.
\end{thm}

\begin{rem} In fact, the assertion that the canonical basis of $\O_K$ is WR twistable if and only if $K=\que(\sqrt{5})$ also follows from Corollary~2 of~\cite{taoufiq}.
\end{rem}

The set of similarity classes of planar WR lattices is parameterized by a curve, while the set of similarity classes of planar stable lattices is two-dimensional (see~\cite{lenny:florian} for details). The results of~\cite{lenny:florian} also indicate that the set of stable ideal lattices in the plane is likely order of magnitude larger than the set of WR ideal lattices: at least this is true for planar arithmetic lattices. With this in mind, it is somewhat surprising that the number of ideal lattices with stable twistable canonical basis is not much larger than with WR twistable canonical basis: this is an interesting property of the canonical basis.

\begin{thm} \label{main_stable} Let $K = \que(\sqrt{D})$ be a real quadratic number field and $\O_K$ its ring of integers. There can be at most finitely many ideals in $\O_K$, up to similarity of the resulting ideal lattices, with the canonical basis being stable twistable: if this is the case for some ideal $I$ with canonical basis $a, b+g\delta$, then
$$b < \left\{ \begin{array}{ll}
\frac{2}{\sqrt{3}} g \sqrt{D} & \mbox{if $D \not\equiv 1\ (\md 4)$,} \\
\frac{g}{2} \left( \frac{2}{\sqrt{3}} \sqrt{D} - 1 \right) & \mbox{if $D \equiv 1\ (\md 4)$.}
\end{array}
\right.$$
The canonical basis for $\O_K$ is stable twistable if and only if $K=\que(\sqrt{5})$.
\end{thm}

\noindent
While our results prove finiteness of the number of ideals in each fixed number field with WR or stable twistable canonical basis, there are still more ideals with stable twistable than with WR twistable canonical basis as we demonstrate in Section~\ref{stable_proof}. We set some additional notation in Section~\ref{setup}, including convenient explicit criteria to check if the canonical basis of a planar ideal lattice is WR or stable twistable, or not. We prove Theorem~\ref{main_wr} and give some explicit examples of ideal lattices with WR twistable canonical basis in Section~\ref{wr_proof}. In Section~\ref{stable_proof} we prove Theorem~\ref{main_stable} and show examples of ideal lattices with stable twistable (but not WR twistable) canonical basis in Section~\ref{wr_proof}. Finally, in Section~\ref{remarks} we make some remarks on the relation between stable and WR twists of arbitrary ideal bases. 

Let us conclude this section with a few additional words of motivation for our results. Our theorems reveal certain new properties of the canonical basis for ideals in quadratic number fields. Canonical bases have an advantage of being convenient for computations, including explicit computations of WR and stable lattices, as demonstrated in~\cite{kate-me}, \cite{fletcher-jones} and~\cite{lenny:stable}. Knowing which of those that are not WR or stable to start with can be twisted into WR or stable shortest bases is helpful. Further, this information may find some future applications when choosing ideal lattice bases for lattice codes constructions. In Section~\ref{comm_theory} we include some additional details on applications of the WR and stable twists of lattice bases in communication theory and in arithmetic theory of quadratic forms and lattices. We are now ready to proceed.
\bigskip

\section{Notation and setup}
\label{setup}

We start here with a basic review of some general lattice theory background. If rank of $L$ is $n \geq 2$, we define the {\it successive minima} of $L$ to be the real numbers
$$0 < \lambda_1 \leq \dots \leq \lambda_n$$
such that
$$\lambda_i = \inf \left\{ \mu \in \real : \dim_{\real} \{ \bx \in L : \|\bx\| \leq \mu \} = i \right\}.$$
Then $\lambda_1 = |L|_{\|\ \|^2}$, and $L$ is WR if and only if $\lambda_1 = \dots = \lambda_n$. Linearly independent vectors $\bx_1,\dots,\bx_n \in L$ such that $\|\bx_i\| = \lambda_i$ are called {\it vectors corresponding to successive minima}. They do not necessarily form a basis for $L$. On the other hand, $L$ has a {\it Minkowski reduced basis} $\bv_1,\dots,\bv_n$, defined by the conditions that
$$\|\bv_1\| = \lambda_1,\ \|\bv_i\| = \min \left\{ \|\bx\| : \bv_1,\dots\bv_{i-1}, \bx \text{ are extendable to a basis for } L \right\}.$$
In general, $\|\bv_i\| \geq \lambda_i$, and when $n \geq 5$ these inequalities can be strict, although a theorem of van der Waerden asserts that for all $i \geq 4$,
\begin{equation}
\label{vdw}
\|\bv_i\| \leq \left( \frac{5}{4} \right)^{i-4} \lambda_i,
\end{equation}
and there is a conjecture that the constant in the upper bound of~\eqref{vdw} can be improved. When $n \leq 4$, a Minkowski reduced basis for a lattice $L$ always consists of vectors corresponding to successive minima. For each $n \geq 2$, there are finitely many inequalities that have to be satisfied by the vectors $\bv_1,\dots,\bv_n \in L$ to be a Minkowski reduced basis for $L$. The number of such inequalities depends only on $n$, but it grows fast with $n$. The explicit list of non-redundant inequalities is known only for $n \leq 7$ (it is called Tammela's list); see \S\S 2.2-2.3 of \cite{achill_book} for Tammela's list, as well as more information on Minkowski reduction and relation to successive minima.

Let $L = (\bv_1\ \bv_2)\ \zed^2$ be a planar lattice. The necessary and sufficient conditions for the basis $\bv_1,\bv_2$ to be Minkowski reduced (and hence to correspond to successive minima $\lambda_1,\lambda_2$) are as follows:
\begin{equation}
\label{mink-2}
\|\bv_1\| \leq \|\bv_2\|,\ 2 | \bv_1^t \bv_2 | \leq \|\bv_1\| \|\bv_2\|. 
\end{equation}
In order for $L$ to be WR we need $\lambda_1=\lambda_2$, i.e.
\begin{equation}
\label{mink-2-WR}
\|\bv_1\| = \|\bv_2\|,\ 2 | \bv_1^t \bv_2 | \leq \|\bv_1\|^2, 
\end{equation}
and in order for $L$ to be stable we need 
\begin{equation}
\label{mink-2-stable}
\sqrt{\det(L)} \leq \|\bv_1\|,\ 2 | \bv_1^t \bv_2 | \leq \|\bv_1\| \|\bv_2\|.
\end{equation}

The angle $\theta$ between $\bv_1$ and $\bv_2$ must therefore be in the interval $[\pi/3,2\pi/3]$, and $| \cos \theta |$ is an invariant of the lattice. Two planar WR lattices $L_1,L_2$ are similar if and only if these values are the same (see~\cite{hex}). There is an easy way to ``deform" a non-WR lattice into a WR one.

\begin{lem} \label{WR_stretch} Let $L = (\bv_1\ \bv_2)\ \zed^2 \subset \real^2$ be a lattice of full rank, where $\bv_1, \bv_2$ is a Minkowski reduced basis matrix, so $\|\bv_1\| = \lambda_1, \|\bv_2\| = \lambda_2$ for the successive minima $\lambda_1 \leq \lambda_2$ of $L$. Let $\bu_1 = \lambda_2 \bv_1$, $\bu_2 = \lambda_1 \bv_2$, then the lattice $M = ( \bu_1\ \bu_2)\ \zed^2$ is a well-rounded lattice with successive minima equal to $\lambda_1 \lambda_2$.
\end{lem}

\proof
This follows immediately from Lemma~3.6 of~\cite{hex}.
\endproof

\begin{rem} In principle, there is a deformation of a lattice into a WR lattice in higher dimensions too, but it is more complicated. Such a deformation is described in Remark~3.3 of~\cite{lizhen_ji}.
\end{rem}
\smallskip

We now come back to the ideal lattices. Let $K=\que(\sqrt{D})$ be a real quadratic number field and $\alpha = p + q \sqrt{D} \in K$ be totally positive, then
$$A(\alpha) = \begin{pmatrix} \sqrt{p+q \sqrt{D}} & 0 \\ 0 & \sqrt{p - q\sqrt{D}} \end{pmatrix},$$
where $p \pm q\sqrt{D} > 0$; in fact, we can assume without loss of generality that $p,q$ are positive rational numbers. Let $B$ as in~\eqref{B1} or~\eqref{B2} be the canonical basis matrix for an ideal lattice $L_K(I)$, where $I \subseteq \O_K$ is the corresponding ideal. Then $B$ is WR twistable (respectively, stable twistable) if 
$$C := A(\alpha) B = \left\{ \begin{array}{ll}
\begin{pmatrix} a \sqrt{p+\sqrt{D}} & (b - g\sqrt{D}) \sqrt{p+\sqrt{D}} \\ a \sqrt{p-\sqrt{D}} & (b + g\sqrt{D}) \sqrt{p-\sqrt{D}} \end{pmatrix} & \mbox{if $D \not\equiv 1\ (\md 4)$,} \\
\begin{pmatrix} a \sqrt{p+\sqrt{D}} & \left( \frac{(2b+g) - g\sqrt{D}}{2} \right) \sqrt{p+\sqrt{D}} \\ a \sqrt{p-\sqrt{D}} & \left( \frac{(2b+g) + g\sqrt{D}}{2} \right) \sqrt{p-\sqrt{D}} \end{pmatrix} & \mbox{if $D \equiv 1\ (\md 4)$.}
\end{array}
\right.$$
is a Minkowski reduced basis for WR (respectively, stable) lattice $C\zed^2$. These conditions can be described by explicit inequalities, stemming from~\eqref{mink-2-WR} and~\eqref{mink-2-stable}. First assume that $D \not\equiv 1 (\md 4)$, then $B$ is WR twistable by $\alpha = p + q\sqrt{D}$ if and only if
\begin{eqnarray}
\label{C1}
a^2p = (Dg^2 + b^2)p - 2qbgD,\ \left| b^2-g^2D+a^2 \right| \leq ab,
\end{eqnarray}
and $B$ is stable twistable by $\alpha$ if and only if
\begin{eqnarray}
\label{C2}
-4 D ^2 g^2 q^2 + 6 b D g p q + D g^2 p^2 - 3 b^2 p^2 \geq 0, \nonumber \\
\min \left\{ a^2p, (Dg^2+b^2)p-2Dbgq \right\} \geq ag \sqrt{D(p^2 - q^2D)}.
\end{eqnarray}
On the other hand, if $D \equiv 1 (\md 4)$, then $B$ is WR twistable by $\alpha = p + q\sqrt{D}$ if and only if
\begin{eqnarray}
\label{Cd1}
2a^2p = \frac{1}{2} \left( (D+1)p - 2Dq) \right) g^2 - 2b(Dq-p) g+2b^2 p, \nonumber \\
\left| 4a^2 + 4b^2 + 4bg - (D-1) g^2 \right| \leq 2a(2b+g),
\end{eqnarray}
and $B$ is stable twistable by $\alpha$ if and only if
\begin{eqnarray}
\label{Cd2}
-D^2 g^2 q^2 + 3 D g b p q + \frac{3}{2} D g^2 p q - 3 b^2 p^2 - 3 b g p^2 - \frac{3}{4} g^2 p^2 + \frac{1}{4} D g^2 p^2 \geq 0 \nonumber \\
\min \left\{ 2 a^2p, ag (p - Dq) + 2abp \right\} \geq ag \sqrt{D(p^2 - q^2D)}.
\end{eqnarray}
We can now use these criteria to analyze the WR and stable twistable properties of the canonical bases for ideals in real quadratic number fields.
\bigskip

\section{Proof of Theorem~\ref{main_wr}}
\label{wr_proof}

In this section we prove Theorem~\ref{main_wr} step by step. Let $K=\que(\sqrt{D})$ for a squarefree integer $D > 1$. Let $I \subseteq \O_K$ be an ideal with the canonical basis $a, b + g\delta$ as described in~\eqref{abg}. Suppose this basis is WR twistable by some totally positive $\alpha = p + q \sqrt{D} \in K$. First assume that $K = \que(\sqrt{D})$ with $D \not\equiv 1\ (\md 4)$, then by~\eqref{C1} we must have
$$p = \left( \frac{2bgD}{b^2+g^2D-a^2} \right) q \text{ and } \left| b^2-g^2D+a^2 \right| \leq ab.$$
We should remark that the first identity holds, unless the denominator $b^2+g^2D-a^2 = 0$, which is not possible: if this is the case, then~\eqref{C1} implies that $b=0$, and so $a^2=g^2D$, which contradicts $D$ being squarefree. Since $p$ must be positive, we have $a^2 < b^2+g^2D$. If $b^2 > g^2D$, we have
$$a^2 < (b^2 - g^2D)  + a^2 \leq ab,$$
meaning that $a < b$, which contradicts the choice of the canonical basis. Hence we must have $b^2 < g^2D$: equality is not possible since $D$ is squarefree. If $I = \O_K$, then $a=1, b=0, g=1$, and so we must have $|1-D| \leq 0$, which is a contradiction. Hence $\O_K$ cannot have WR twistable canonical basis.

Now suppose that $D \equiv 1\ (\md 4)$, then by~\eqref{Cd1} we must have
$$p = \left( \frac{2gD (2b + g)}{4b^2+g^2(D+1)+4b-4a^2} \right) q,$$
unless $4b^2+g^2(D+1)+4b-4a^2 = 0$, in which case $2b + g=0$: this is not possible, since $g \mid b$. Additionally,~\eqref{Cd1} implies that
\begin{equation}
\label{Cd-3}
\left| 4a^2 + 4b^2 + 4bg - (D-1)g^2 \right| \leq 2a (2b+g).
\end{equation}
Since $p$ must be positive, we have $a^2 < b^2+b+\frac{g^2(D+1)}{4}$.  If $b^2 + bg > \frac{(D-1)g^2}{4}$, then we get
$$4a^2 < 4a^2 + 4b^2 + 4bg - (D-1)g^2 \leq 2a (2b+g),$$
and so $a < b+g/2$, which contradicts the fact that $g \mid a-b$. Hence we must have $b^2 + bg \leq \frac{(D-1)g^2}{4}$, which means that $b < \frac{(\sqrt{D}-1)g}{2}$: again, equality is not possible since $D$ is squarefree. If $I = \O_K$, then $a=1, b=0, g=1$, so $p = \frac{2Dq}{D-3}$. Hence $\alpha = q \left( \frac{2D}{D-3} + \sqrt{D} \right)$, which again is not totally positive unless $D=5$ (otherwise $D \geq 13$, and so $2D/(D-3) < \sqrt{D}$). If $D=5$, then we can take $\alpha = 5 + \sqrt{5}$, obtaining the matrix
$$A(\alpha) B = \begin{pmatrix} \sqrt{5+\sqrt{5}} & \frac{(1-\sqrt{5}) \sqrt{5+\sqrt{5}}}{2} \\ \sqrt{5-\sqrt{5}} & \frac{(1+\sqrt{5}) \sqrt{5-\sqrt{5}}}{2} \end{pmatrix}$$
with orthogonal columns, both of norm $= \sqrt{10}$. Hence $\O_K$ cannot have WR twistable canonical basis for any real quadratic number field $K \neq \que(\sqrt{5})$. 

Finally notice that if $I$ is an ideal with canonical basis $a, b+g\delta$ and $I' = \frac{1}{g} I$ is the corresponding ideal with canonical basis $\frac{a}{g}, \frac{b}{g} + \delta$, then the lattices $L_K(I)$ and $L_K(I')$ are similar, and so are there twists by the same element~$\alpha \in K$. Therefore our upper bounds on $b$ mean that there can be at most finitely many ideals in $\O_K$, up to similarity of the resulting ideal lattices, with the canonical basis being WR twistable.
\medskip

This finishes the proof of Theorem~\ref{main_wr}. Notice that this theorem implies that the canonical basis is rarely WR twistable, however such examples still exist. We demonstrate a couple here.

\begin{ex} \label{ex1} Let $D = 139 \equiv 3\ (\md 4)$ and $K=\que(\sqrt{139})$. Let $I \subset \O_K$ be an ideal generated by the canonical basis
$$9,\ 7-\sqrt{139},$$
i.e. $g=1$, $b=7 < \sqrt{139}$, $a = 9 \mid \Nn(b-\sqrt{D}) = 7^2-139 = -90$. This basis is WR twistable by the totally positive element
$$\alpha = \frac{1946}{107} + \sqrt{139} \in K$$
with the resulting WR lattice having the Minkowski reduced basis matrix
$$\frac{1}{107} \begin{pmatrix} 9 \sqrt{208222+11449 \sqrt{139}} & (7 - \sqrt{139}) \sqrt{208222+11449 \sqrt{139}} \\ 9 \sqrt{208222-11449 \sqrt{139}} & (7+\sqrt{139}) \sqrt{208222-11449 \sqrt{139}}  \end{pmatrix}$$
and cosine of the angle between these basis vectors being $-1/14$. The common value of the successive minima of this lattice is $\sqrt{\frac{315252}{107}} \approx 54.27964973...$
\end{ex}

\begin{ex} \label{ex2} Let $D = 141 \equiv 1\ (\md 4)$ and $K=\que(\sqrt{141})$. Let $I \subset \O_K$ be an ideal generated by the canonical basis
$$5,\ 4+\frac{1-\sqrt{141}}{2},$$
i.e. $g=1$, $b=4 < \sqrt{141}/2$, $a = 5 \mid \Nn\left( b+\frac{1-\sqrt{D}}{2} \right) = \frac{(8+1)^2-141}{4} = -15$. This basis is WR twistable by the totally positive element
$$\alpha = \frac{1269}{61} + \sqrt{141} \in K$$
with the resulting WR lattice having the Minkowski reduced basis matrix
$$\frac{1}{61} \begin{pmatrix} 5 \sqrt{77409+3721 \sqrt{141}} & \left( \frac{9 - \sqrt{141}}{2} \right) \sqrt{77409+3721 \sqrt{141}} \\ 5 \sqrt{77409-3721 \sqrt{141}} & \left( \frac{9 + \sqrt{141}}{2} \right) \sqrt{77409-3721 \sqrt{141}}  \end{pmatrix}$$
and cosine of the angle between these basis vectors being $2/9$. The common value of the successive minima of this lattice is $\sqrt{\frac{63450}{61}} \approx 32.25157258...$
\end{ex}
\bigskip

\section{Proof of Theorem~\ref{main_stable}}
\label{stable_proof}

Here we prove Theorem~\ref{main_stable}. As above, let $K=\que(\sqrt{D})$ for a squarefree integer $D > 1$ and let $I \subseteq \O_K$ be an ideal with the canonical basis $a, b + g\delta$ as described in~\eqref{abg}. Suppose this basis is stable twistable by some totally positive $\alpha = p + q \sqrt{D} \in K$. First assume that $K = \que(\sqrt{D})$ with $D \not\equiv 1\ (\md 4)$. If we consider the first condition of~\eqref{C2} as a quadratic inequality in $p$, then its leading coefficient is negative unless $b \leq g \sqrt{D/3}$ and it has negative discriminant unless $b \leq 2g\sqrt{D/3}$. Hence the canonical basis can be stable twistable only if 
$$b < \frac{2}{\sqrt{3}} g \sqrt{D}.$$
As always, equality is not possible since $b$ is an integer and $D$ is squarefree. Now assume that $I=\O_K$, then $a=1,b=0,g=1$, so the inequalities of~\eqref{C2} become:
$$D(p^2-4Dq^2) \geq 0,\ p \geq \sqrt{D(p^2-Dq^2)}.$$
Combining these inequalities, we obtain $(1-3D)p^2 \geq 0$, which is not possible. Hence $\O_K$ cannot have a stable twistable canonical basis.
\medskip

Now suppose that~$D \equiv 1\ (\md 4)$. Considering the first condition of~\eqref{Cd2} as a quadratic inequality in $p$, we see that its leading coefficient is negative unless $b \leq \frac{g}{2} \left( \sqrt{D/3} - 1 \right)$ and it has negative discriminant unless $b \leq \frac{g}{2} \left( 2 \sqrt{D/3} - 1 \right)$. Hence the canonical basis can be stable twistable only if 
$$b < \frac{g}{2} \left( 2 \sqrt{D/3} - 1 \right).$$
Again, equality is not possible since $b$ is an integer and $D$ is squarefree. Now assume that $I=\O_K$, then $a=1,b=0,g=1$, so the inequalities of~\eqref{Cd2} become:
$$\frac{p^2}{4} (D-3) + \frac{3}{2} D pq - D^2 q^2 \geq 0,\ 2p > \sqrt{D (p^2-q^2D)}.$$
These inequalities lead to 
$$q \left( \frac{D \sqrt{4D-3}-3}{D-3} \right) < p < q \left( \frac{D}{\sqrt{D-4}} \right),$$
which means that we must have $\frac{D \sqrt{4D-3}-3}{D-3} < \frac{D}{\sqrt{D-4}}$. This only holds for $D=5$, which, as we know, yields even a WR twist. By the same reasoning as in the proof of Theorem~\ref{main_wr}, our upper bounds on $b$ mean that there can be at most finitely many ideals in $\O_K$, up to similarity of the resulting ideal lattices, with the canonical basis being stable twistable. This finishes the proof of Theorem~\ref{main_stable}. 
\bigskip

Finally, let us demonstrate a couple examples of stable twists of ideal lattice canonical bases that are not WR twistable.

\begin{ex} \label{ex1-1} Let $D = 1327 \equiv 3\ (\md 4)$ and $K=\que(\sqrt{1327})$. Let $I \subset \O_K$ be an ideal generated by the canonical basis
$$39,\ 38-\sqrt{1327},$$
i.e. $g=1$, $b=38$, $a = 39 \mid \Nn(b-\sqrt{D}) = 38^2-1327 = 117$. Notice, in particular, that
$$g \sqrt{D} = 36.42801120... < b < 42.06344416... = \frac{2}{\sqrt{3}} g \sqrt{D},$$
hence this basis cannot be WR twistable, by Theorem~\ref{main_wr}. On the other hand, this basis is stable twistable by the totally positive element
$$\alpha = 63 + \sqrt{1327} \in K$$
with the resulting stable lattice having the Minkowski reduced basis matrix
$$\begin{pmatrix} 39 \sqrt{63 + \sqrt{1327}} & (38 - \sqrt{1327}) \sqrt{63 + \sqrt{1327}} \\  39 \sqrt{63 - \sqrt{1327}} & (38 + \sqrt{1327}) \sqrt{63 - \sqrt{1327}} \end{pmatrix}$$
and cosine of the angle between these basis vectors being $0.4951063950...$. The determinant of this lattice is $146048.2881...$ and values of the successive minima are $\sqrt{147442}$ and $\sqrt{191646}$, respectively.
\end{ex}

\begin{ex} \label{ex1-2} Let $D = 125173 \equiv 1\ (\md 4)$ and $K=\que(\sqrt{125173})$. Let $I \subset \O_K$ be an ideal generated by the canonical basis
$$183,\ 182 + \frac{1-\sqrt{125173}}{2},$$
i.e. $g=1$, $b=182$, $a = 183 \mid \Nn \left( \frac{(2b+1) - \sqrt{D}}{2} \right) = 2013$. Notice, in particular, that
$$\frac{g (\sqrt{D}-1)}{2} = 176.3989824... < b < 203.7653503... = \frac{g}{2} \left( \frac{2}{\sqrt{3}} \sqrt{D} - 1 \right),$$
hence this basis cannot be WR twistable, by Theorem~\ref{main_wr}. On the other hand, this basis is stable twistable by the totally positive element
$$\alpha = 611 + \sqrt{125173} \in K$$
with the resulting stable lattice having the Minkowski reduced basis matrix
$$\begin{pmatrix} 183 \sqrt{611 + \sqrt{125173}} & \left( \frac{365}{2} - \frac{\sqrt{125173}}{2} \right) \sqrt{611 + \sqrt{125173}}  \\  183 \sqrt{611 - \sqrt{125173}} & \left( \frac{365}{2} + \frac{\sqrt{125173}}{2} \right) \sqrt{611 - \sqrt{125173}} \end{pmatrix}$$
and cosine of the angle between these basis vectors being $0.4853755919...$. The determinant of this lattice is $32252383.1...$ and values of the successive minima are $\sqrt{33252444}$ and $\sqrt{40923558}$, respectively.
\end{ex}
\medskip

\section{Further remarks on stable and WR twistable bases }
\label{remarks}

\begin{figure}[H]
\centering
\includegraphics[scale=0.4]{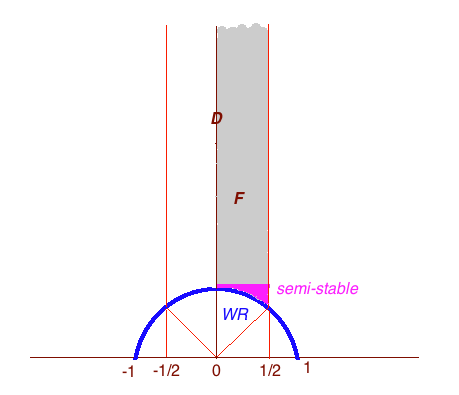}
\caption{Similarity classes of lattices in $\real^2$ with WR and stable subregions marked by colors.}\label{fig:domain}
\end{figure}

Here we include some further heuristic remarks on stable and WR twistable bases beyond the canonical ones discussed above. Let
$$\hyp = \{ x+iy \in \cee : y>0 \}$$
be the upper half-plane. Following~\cite{lenny:florian}, define
$$\D := \{ \tau = a + bi \in \hyp : -1/2 < a \leq 1/2, |\tau| \geq 1 \}$$
and
$$\F := \{ \tau = a + bi \in \hyp : 0 \leq a \leq 1/2, |\tau| \geq 1 \},$$
so $\F$ is ``half" of $\D$. Then $\D$ is the standard fundamental domain of $\hyp$ under the action of $\SL_2(\zed)$ by fractional linear transformations and $\F$ is the space of similarity classes of planar lattices. As in~\cite{lenny:florian}, this can be illustrated by Figure~\ref{fig:domain}  with the subsets of WR and stable similarity classes marked accordingly.

Let $K$ be a real quadratic field and $I \subseteq \O_K$ an ideal. For each lattice $(L_K(I),f_{\A(\alpha)})$ for a totally positive $\alpha \in K$, let us write $\left< L_K(I),f_{\A(\alpha)} \right>$ for its similarity class. As shown in~\cite{bayer-nebe}, sets of similarity classes of the form
$$\G(I):= \left\{ \left< L_K(I),f_{\A(\alpha)} \right> : \alpha \in K \textrm{ totally positive} \right\}$$
correspond to closed geodesics in $\hyp/\SL_2(\zed)$, and every closed geodesic will necessary intersect the WR locus $\W$ in the upper-half plane (the blue arc in Figure~\ref{fig:domain}); denote by $\I := \W \cap \G(I)$ the set of such intersection points. 

\begin{figure}[H]
\centering
\includegraphics[scale=0.4]{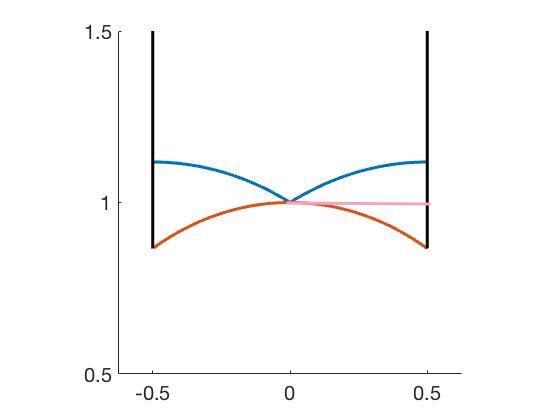}
\caption{Geodesic $\G(I)$ (blue) intersects the WR locus (red) in one point (orthogonal lattice), which is also a stable twist.}\label{WR-orth}
\end{figure}

The set $\I$ is completely defined by a relation described in~\cite{taoufiq}. Let $B=\{x,y\}$ be a basis for the ideal $I$, and define 
$$F(B) = \Nn(x)^2+\Nn(y)^2 + \Nn(x) \Nn(y)- \Nn(I)^2 \Delta_K /4,$$
where $\Delta_K$ is the discriminant of the number field $K$. We say that two bases $B$ and $B'$ are {\it equivalent} if $F(B)= F(B')$: this is indeed an equivalence relation. Then the set $\I$ is in bijective correspondence with the equivalence classes of bases $B$ with $F(B) < 0$. This implies  finiteness of $|\I|$. Furthermore, as the geodesics $\G(I)$ are closed, the number of arcs ``inside" the stable locus is at most twice $|\mathcal{I}|$, except when the only WR twist of our ideal lattice is orthogonal as in Figure~\ref{WR-orth}.

For example, if $I=\O_K$ this will happen if and only if $D=s^2 +1$ (respectively $D=s^2 +4$) for $D \not\equiv 1\pmod 4$ (respectively $D \equiv 1 \pmod 4$). In general, this will be the case if for every basis $B=\{x,y\}$ of $I$ such that $\Nn(x) < \Nn(I) \sqrt{\Delta_K}/4$, we have $\Nn(x) = -\Nn(y)$.

\begin{figure}[H]
\centering
\includegraphics[scale=0.4]{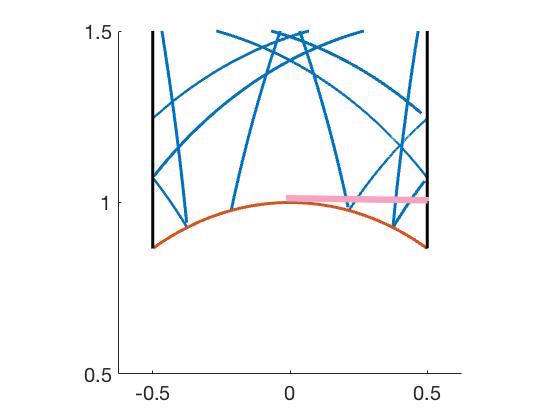}
\caption{Geodesic $\G(\O_{\que(\sqrt{59})})$ intersects the WR locus.}\label{WR-59}
\end{figure}

Figure~\ref{WR-59} illustrates an example of the intersection of $\G(I)$ with the WR locus for $I=\O_K$, where $K=\que(\sqrt{59})$.

Before intersecting the WR locus, the geodesic $\G(I)$ crosses ``continuously" the stable locus, which yields infinitely many stable twists: in other words, this insures the existence of infinitely many totally positive $\alpha \in K$ such that $\A(\alpha)B$ is stable for a fixed WR twistable basis.  
In summary, every WR twistable basis will give rise to infinitely many stable twists, except in the cases where the ideal $I$ admits only the orthogonal WR twist. The converse is not true: there exist bases that are stable twistable, but not WR twistable, as demonstrated by Examples~\ref{ex1-1} and~\ref{ex1-2}.
\bigskip

\section{Applications of WR and stable lattice twists}
\label{comm_theory}

In this section we discuss some applications of WR and stable twistable lattice bases. First we highlight a connection between such bases and the theory of error control in wireless communication. We assume a single-input single-output (SISO) Rayleigh flat fading channel model:
\begin{equation}
\label{model}
\bwy = H\bx + \bv,  
\end{equation}
where $\bx$ is the transmitted codeword taken from some finite codebook in $\cee^n$, $H = \diag(h_1,\dots, h_n)$ describes the random channel response, and $\bv \in \cee^n$ is a random additive white Gaussian noise with variance $\sigma_{\bv}^2$.

To study the communication reliability of a given code $C$ we consider the codeword error probability $P_e(C)$. The goal is to choose $C$ to be a subset of a lattice that minimizes $P_e(C)$. Considering $\Lambda \subset \real^n$, it is proved in \cite{oggier} that the lattices (of fixed volume) that minimize $P_e(C)$ for Rayleigh fading channel are those with maximal $d_{\min}(\Lambda)$, where
$$d_{\min}(\Lambda) = \min_{\bx \in \Lambda \setminus \{\bo\}} \prod_{i=1}^n |x_i|.$$
This criterion is restricted to the so called fully diverse lattices, i.e. the lattices with non-vanishing~$d_{\min}(\Lambda)$. Let us prove that the existence of twistable bases allows one to restrict this optimization problem to the set of WR or stable lattices without loss of generality.

\begin{prop}
Let $\Lambda$ be a lattice with maximal $d_{\min}(\Lambda)$ in its dimension. Then there exists a WR lattice $L$ and a stable lattice $M$, such that 
$$d_{\min}(\Lambda)=d_{\min}(L)=d_{\min}(M).$$
\end{prop}

\proof
Let $\Lambda\subset \real^n$ be a fully diverse lattice. Notice that $d_{\min}(\Lambda)$ is invariant under the action of $\A$, i.e $d_{\min}(A \Lambda) = d_{\min}(\Lambda)$ for any $A \in \A$. In \cite{mcmullen} Mcmullen showed that if the orbit closure $\overline{\A\Lambda}$ is compact then $\A\Lambda$ meets the set of WR lattices. Hence we only need to show that the full diversity of $\Lambda$ will ensure the compactness of $\overline{\A.\Lambda}$, and this is a straightforward application of the Mahler compactness criterion. Namely, for a set $E$ of unimodular lattices we have
$$\overline{E}\textrm{ is compact} \iff \lambda_1(L)>0 \textrm{ for all } L\in E.$$
If $\Lambda$ is a fully diverse lattice, then $d_{\min}(\Lambda) > 0$. Hence, by the AM-GM inequality, we have  
$$0 < \sqrt{n} \prod_{i=1}^{n} |x_i| \leq \|\bx\|$$
for any $\bx \in \Lambda$. In particular, taking $\bx \in \Lambda$ such that $\|\bx\| = \lambda_1(\Lambda)$, we see that $\lambda_1(L) > 0$ for any $L\in \A \Lambda$.

For stable lattices the argument is the same replacing the result of~\cite{mcmullen} with the analogous result of~\cite{weiss} which guarantees that if the orbit closure $\overline{\A \Lambda}$ is compact then $\A\Lambda$ meets the set stable lattices.
\endproof

\begin{rem} Notice that $d_{\min}(\Lambda)>0$ for $\Lambda$ arising from a real number field. Margulis conjectured that the converse is also true. More precisely, if a lattice $\Lambda \subset \real^n$ for $n \geq 3$ has $d_{\min}(\Lambda) >0$ then $\Lambda$ comes from an order in a number field.
In other words, Margulis's conjecture asserts that all fully diverse lattices come from orders in totally real number fields, and this motivates the choice of number theoretic constructions in this context.
\end{rem}
\medskip

Let us also briefly discuss the use of WR and stable twists in the arithmetic theory. Let $K$ be a real number field of degree~$n$. The Euclidean minimum of $K$ is defined as
$$M(K) := \inf \left\{ \alpha \in \real_{>0} : \forall x \in K~\exists y \in \O_K \textrm{ such that } |\Nn(x-y)| \leq \alpha \right\}.$$
The Euclidean minimum measures how far is $K$ from having a Euclidean algorithm. In fact, if $M(K) < 1$ then $\O_K$ is a Euclidean ring. Furthermore, for an ideal $I$ in $\O_K$ we can define 
$$M(I) := \inf \left\{ \alpha \in \real_{>0} : \forall x \in K~\exists y \in I \textrm{ such that } |\Nn(x-y)| \leq \alpha \right\}.$$
Now, for a lattice $\Lambda \subset \real^n$, its covering radius is 
$$\mu(\Lambda) := \sup_{\bwy \in \real^n} \min \left\{ \| \bx - \bwy \| : \bx \in \Lambda \right\}$$
and the Hermite thickness of $\Lambda$ is 
$$\tau(\Lambda) := \frac{\mu(\Lambda)^2}{\det(\Lambda)^{1/n}}.$$
Then for an ideal $I \subseteq \O_K$, define
$$\tau_{\min}(I):= \min\{\tau(L_K (I),f_{A(\alpha)}) : \alpha \in K~\textrm{totally positive}\}.$$

\begin{thm} [\cite{bayer}]
For all number fields $K$ of degree $n$, 
$$M(I) \leq \left( \frac{\tau_{\min}(I)}{n} \right)^{n/2} \sqrt{|\Delta_K|}\ \Nn(I),$$
where $\Delta_K$ is the discriminant of $K$.
\end{thm}

This theorem motivates the study of the covering radii in the orbit of the action of $\A_1(K)$ on $L_K (I)$.
Furthermore, it has been proved that $\mu(\Lambda) \leq \frac{\sqrt{n}}{2}$ for any unimodular WR lattice of dimension $n\leq 9$ (not true for $n\geq 30$), and it is conjectured to be true for stable lattices in any dimension.
Combining these observations, we have
$$M(K) \leq \frac{\sqrt{|\Delta_K|}}{4}$$
when $K$ is a real quadratic number field. Notice that this bound automatically implies that $\mathbb{Q}(\sqrt{d})$ is a Euclidean domain for $d=5,2,3,13$.
\medskip

{\bf Acknowledgment.} We would like to thank Laia Amoros, Camilla Hollanti and David Karpuk for many helpful conversations on the subject of this paper. A part of the work was done during a stay of the first author at Claremont McKenna College (CMC): M.T. Damir would like to thank CMC for its hospitality. We also thank the anonymous referee for some helpful suggestions.

\bibliographystyle{plain}  

\end{document}